%% file: lip.tex
\input rbmacros.tex
\def\sS{{\cal S}}

\def\bone{{\bf 1}}
\def\df{{\mathop {\ =\ }\limits^{\rm{df}}}}

\centerline{\bf THE SUPREMUM OF BROWNIAN LOCAL TIMES}
\cl{\bf ON H\"OLDER CURVES}
\footnote{$\empty$}{\rm Research partially supported by NSF grant DMS-9700721.}

\vskip0.4truein
\centerline{{\bf Richard F.~Bass}
\qq and \qq
{\bf Krzysztof Burdzy}}

\vskip1truein

{\narrower
{\bf Abstract}. For $f: [0,1]\to \R$, we consider $L^f_t$,
the local time of space-time Brownian motion on the curve
$f$. 
Let $\sS_\al$ be the class of all functions whose H\"older norm of order $\al$
is less than or equal to 1. We show that the supremum of $L^f_1$ over
$f$ in $\sS_\al$ is finite is $\al>\frac12$ and infinite if $\al<\frac12$.
}

\vfill\eject

\subsec{1. Introduction}

Let $W_t$ be one-dimensional Brownian motion and let $f:[0,1]\to \R$
be a H\"older continuous function. There are a number of equivalent ways to define the 
local time  of $W_t$ along the curve $f$. We will show the equivalence
below, but for now define $L^f_t$ as the limit in probability of
$$\frac{1}{2\eps}\int_0^t \bone_{(f(s)-\eps, f(s)+\eps)}(W_s)\, ds$$
as $\eps\to 0$. Let
$$\sS_\al=\{f: \sup_{0\leq t\leq 1}|f(t)|\leq 1,  
|f(s)-f(t)|\leq |s-t|^\al\hbox{ if }s,t\leq 1\}.$$

We were led to the results in this paper  by the following question.

\proclaim Question 1.1. Is $\sup_{f\in \sS_1} L^f_1$ finite or infinite?

\bigskip

Our interest in this problem arose when we were working on Bass and Burdzy (1999).
A positive answer to Question 1.1 at that time would have provided a proof
of uniqueness for a certain stochastic differential equation; we ended up
using different methods.

However, probably the greatest interest in Question 1.1 has to do with
questions about metric entropy. The metric entropy of $\sS_1$ is known to 
be of order 
$1/\eps$; see, e.g., Clements (1963). That is, if one takes the cardinality of the smallest $\eps$-net
for $\sS_1$ (with respect to the supremum norm) and takes the logarithm,
the resulting number will be bounded above and below by positive constants
times $1/\eps$. It is known (see Ledoux and Talagrand (1991)) that this
is too large for standard chaining arguments to be used to prove finiteness
of $\sup_{f\in \sS_1} L^f_1$.
Nevertheless, the supremum in Question 1.1 is finite.

It is a not uncommon belief among the probability community that metric
entropy estimates are almost always sharp: the supremum of a process is
finite if the metric entropy is small enough, and infinite otherwise.
That is not the case here. 
Informally, our main result is
\bigskip
\noindent{\bf Theorem 1.2}. {\sl The supremum of $f\to L^f_1$
over $\sS_\al$ is finite if $\al >\frac12$ and infinite if $\al < \frac12$.
}
\ms
\ni See Theorems 3.6 and 3.8 for formal statements.

The metric entropy of
$\sS_\al$ when $\al\in (\frac12, 1]$  is far beyond what chaining methods can handle.
Sometimes the method of majorizing measures provides a better result
than that of metric entropy. 
We do not know if this is the case here.

For previous work on local times for space-time curves, see Burdzy and
San Mart\'in (1995) and  Davis (1998). For some results on local times
on Lipschitz curves for two-dimensional Brownian motion, see
Bass and Khoshnevisan (1992) and Marcus and Rosen (1996).

In Section 2 we prove the equivalence of various definitions of $L^f_t$
as well as some lemmas of independent interest. 
In Section
3 we prove finiteness of the supremum over $\sS_\al$ when $\al>\frac12$ and  
that this fails when $\al<\frac12$. We also show that
$(f,t)\to L^f_t$ is jointly continuous on $\sS_\al \times [0,1]$
when $\al > 1/2$.

The letter $c$ with subscripts will denote finite positive constants
whose exact values are unimportant. 
We renumber them in each proof.

\bs
\ni{\bf Acknowledgments} We would like to thank  F.~Gao, E.~Gin\'e, 
J.~Kuelbs, T.~Lyons, and J.~Wellner for their interest and help.
We would like to express our special gratitude to R.~Adler and M.~Barlow
for long discussions of the problem and many instances of specific
advice.

\subsec{2. Preliminaries}

We discuss three possible definitions of $L^f_t$.
\item{(i)} $L^f_t=\lim_{\eps\to 0} \frac{1}{2\eps} 
\int_0^t \bone_{(f(s)-\eps, f(s)+ \eps)}(W_s) ds$;
\item{(ii)} $L^f_t$ is the continuous additive functional of space-time 
Brownian motion associated to the potential $U^f(x,t)=\int_0^{1-t} p(s,x,f(t+s)) ds$,
where $p$ is the transition density for one-dimensional Brownian motion;
\item{(iii)} (for $f\in \sS_1$ only) $L^f_t$ is the local time in
the semimartingale sense at 0 of the process $W_t-f(t)$.

One of the goals of this section is to show the equivalence of these
definitions. We begin with the 
following lemma which will be used repeatedly throughout the paper.

\proclaim Lemma 2.1. Suppose $A^1_t$ and $A^2_t$ are two nondecreasing
continuous processes with $A^1_0=A^2_0=0$. Let $B_t=A_t^1-A_t^2$.
Suppose 
that for all $s\leq t$, and some right-continuous filtration $\{\F_t\}$,
$$\E[A^i_t-A^i_s\mid \F_s]\leq M, \qq \hbox{a.s.}\qq i=1,2,$$
and
for all $s\leq t$
$$\bigl|\,\E[B_t-B_s\mid \F_s]\,\bigr|\leq \gamma, \qq \hbox{a.s.}$$
There exist $c_1,c_2$ such that for all $\lam >0$,
$$\P(\sup_{s\leq t} |B_s|>\lam\sqrt{\gamma M})\leq c_1e^{-c_2\lam}.$$

\proof We have
$$(B_t-B_s)^2= 2\int_s^t (B_t-B_r) dB_r.$$ 
Using  a Riemann sum approximation (cf.~Bass (1995), Exercise I.8.28) we obtain
$$\eqalign{\E&[(B_t-B_s)^2\mid \F_s]
=2\E\Big[\int_s^t (B_t-B_r) \, dB_r\mid \F_s\Big]\cr
&=2\E\Big[\int_s^t \E[B_t-B_r\mid \F_r] \, dB_r\mid \F_s\Big]\cr
&\leq 2\E\Big[\int_s^t \gamma (dA_r^1+dA_r^2)\mid \F_s\Big]\leq
4\gamma M.\cr}$$
This inequality holds a.s.\ for each $s$. The left hand side is equal
to $$\E[B_t^2\mid \F_s]-2B_s\E[B_t\mid \F_s]+B_s^2$$ and hence is right
continuous. Therefore there is a null set outside of which
$$\E[(B_t-B_s)^2\mid F_s]\leq 4\gamma M$$ for all $s$.
In particular, if $T$ is a stopping time, by Jensen's inequality we
obtain
$$\E[|B_t-B_T|\, \mid \F_T]\leq (\E[ (B_t-B_T)^2\mid \F_T])^\half
\leq (4\gamma M)^\half.$$
Our result now follows by Bass (1995, Theorem I.6.11), and Chebyshev's inequality.
\qed

Let $W_t$ be one-dimensional Brownian motion. 
Define $$p(t,x,y)=(2\pi t)^{-\half} \exp(-|x-y|^2/2t),\eqno (2.1)$$
 the transition
density of one dimensional Brownian motion.
In the rest of the paper, $\F_t$ will denote the (right-continuous)
filtration generated by $W_t$.

For a measurable function $f:[0,1]\to \R$
set $\norm{f}=\sup_{t\leq 1}|f(t)|$.
Let $$D^f_t(\eps)=\frac{1}{2\eps}\int_0^t 
\bone_{(f(s)-\eps, f(s)+\eps)}(W_s)\, ds.
$$

\proclaim Proposition 2.2. For $f$ measurable on $[0,1]$,
there exists a nondecreasing continuous process $L^f_t$ such that
$\E \norm{D^f(\eps)-L^f}^2\to 0$ as $\eps\to 0$.

\proof 
Let $\E^{(x,t)} $ denote the expectation corresponding
to the distribution of Brownian motion starting
from $x$ at time $t$, i.e., satisfying $W_t = x$.
For any $x$ and any $t\leq 1$,
$$\eqalignno{\E^{(x,t)} \frac{1}{2\eps}\int_0^{1-t} 
\bone_{(f(t+s)-\eps,f(t+s)+\eps)}(W_{t+s})\, ds
&=\frac{1}{2\eps}\int_0^{1-t} \int_{f(t+s)-\eps}^{f(t+s)+\eps} p(s,x,y)\, dy\, ds
\cr&\leq c_1\int_0^{1-t} \frac{1}{\sqrt{ s}} \, ds\leq c_2\sqrt{1-t}\leq c_2.
&(2.2)\cr}$$
This implies that,
$$\E[ D^f_1(\eps)-D^f_t(\eps)\mid \F_t]=\E^{(W_t,t)}
\frac{1}{2\eps}\int_0^{1-t} \bone_{(f(t+s)-\eps,f(t+s)+\eps)}(W_{t+s})\, ds
\leq c_2.\eqno(2.3)$$ 
The supremum of
$$\frac{1}{2\eps} \int_{f(t+s)-\eps}^{f(t+s)+\eps} p(s,x,y)\, dy
$$ 
over $\eps>0$, $t\leq 1$ and $s\leq 1-t$ is bounded.
By the continuity of $p(s,x,y)$ in $y$ and the bounded convergence
theorem, as $\eps\to 0$,
$$\frac{1}{2\eps}\int_0^{1-t} \int_{f(t+s)-\eps}^{f(t+s)+\eps} p(s,x,y)\, dy\, ds
\to \int_0^{1-t} p(s,x,f(t+s))\, ds$$
uniformly over $x$ and $t$. 
Calculations similar to those in (2.2) and (2.3) yield
the following estimate:  for any $\eta>0$,
$$\left|\E[(D^f_1(\eps_1)-D^f_1(\eps_2)) 
-(D^f_t(\eps_1)-D^f_t(\eps_2))\mid \F_t]\right|
\leq \eta, \qq \hbox{a.s.},\eqno(2.4)$$
for all $t\leq 1$ provided $\eps_1$ and
$\eps_2$ are small enough.

Because of (2.3) and (2.4), we can apply Lemma 2.1 with
$A^1_t = D^f_t(\eps_1)$ and $A^2_t = D^f_t(\eps_2)$.
The estimate in that lemma shows that,
in a sense, the supremum of the difference
between $ D^f_t(\eps_1)$ and $ D^f_t(\eps_2)$
is of order $\sqrt{\eta}$.
We see that $\E(\norm{D^f(\eps_1)-D^f(\eps_2)}^2)\to 0$
as $\eps_1, \eps_2\to 0$. 
This implies that  $\{D^f(\eps_n)\}$ is a Cauchy sequence, and therefore
$D^f(\eps_n)$ converges as $n\to \infty$, for any sequence
$\{\eps_n\}$ converging to $0$. 
Denote the limit by $L^f_t$; it is routine to check that
the limit does not depend on the sequence $\{\eps_n\}$. 
Since the convergence
is uniform over $t$ and $t\to D^f_t(\eps)$
is continuous for every $\eps$, then $L^f_t$ is continuous in $t$.
For a similar reason, $t\to L^f_t$ is nondecreasing.
\qed

\ni{\bf Remark 2.3.}
A very similar proof shows that $L^f_t$ is the limit in
$L^2$ of $$\frac{1}{\eps} \int_0^t \bone_{[f(s), f(s)+\eps)}(W_s) ds.$$

\ms 
\noindent{\bf Remark 2.4}.
Let  
$$U^f(x,t)=\int_0^{1-t}  p(s,x,f(t+s)) \, ds.$$
A straightforward limit argument  shows that
$$\E[L^f_1-L^f_t\mid \F_t]=\int_0^{1-t} p(s,W_t,f(t+s))\, ds.\eqno(2.5)$$
It follows that $U^f(W_t,t)$ is a potential
for the space-time Brownian motion $t\to (W_t,t)$. 
Hence the function $U^f(x,t)$
is excessive with respect to space-time Brownian motion, and therefore
$L^f_t$ can also be viewed as the continuous additive functional
for the space-time Brownian motion $(W_t,t)$ 
whose potential is $U^f$.

\proclaim Corollary 2.5. Suppose $f_n\to f$ uniformly.
Then $ \norm{L^{f_n}-L^f}$ converges to 0 in
$L^2$.

\proof {}From (2.5),
$$\E[L^f_1-L^f_u\mid \F_u]\leq c_1\int_0^{1-u} \frac{1}{\sqrt{ s}}\, ds
\leq c_2 \sqrt{1-u}\leq c_2$$
and
$$\eqalign{\bigl|\, \E[L^{f_n}_1-L^{f_n}_u\mid \F_u]&-\E[L^f_1-L^f_u\mid \F_u]\,\bigr|\cr 
&=\Bigl|\int_0^{1-u} [p(s,W_u,f_n(u+s))-p(s,W_u,f(u+s))]\, ds\Bigr|\cr
&\leq \int_0^{1-u} |p(s,W_u,f_n(u+s))-p(s,W_u,f(u+s))|\, ds.\cr}$$
The right hand side tends to 0 by the assumption that
$f_n\to f$ uniformly, and the result now follows by Lemma 2.1,
using the same argument as at the end of the proof of Proposition 2.2.
\qed

If $f$ is a Lipschitz function, then $W_t-f(t)$ is a semimartingale.
We can therefore define a local time for $W_t$ along the
curve $f$ by setting
$K^f_t$ to be the local time (in the semimartingale sense) at 0 of
$Y_t=W_t-f(t)$. That is,
$$K^f_t=|Y_t|-|Y_0|-\int_0^t \sgn(Y_s)\, dY_s.$$

\proclaim Proposition 2.6. With probability one, $K^f_t=L^f_t$ for all $t$.

\proof By Revuz and Yor (1994) Corollary VI.1.9, 
$$K^f_t=\lim_{\eps\to 0} \frac{1}{\eps}\int_0^t 
\bone_{[0,\eps)}(Y_s) d\angel{Y}_s. \eqno (2.6)$$
Since $Y_t=W_t-f(t)$, then $\angel{Y}_t=\angel{W}_t=t$, and so by
Remark 2.3, $K^f_t=L^f_t$ a.s. Since both $K^f_t$ and $L^f_t$ are continuous
in $t$, the result follows.
\qed

\subsec{3. The supremum of local times}

Our first goal is to obtain an estimate on the number of rectangles
of size $(1/N)\times (2/\sqrt{N})$ that are hit by a Brownian path.
Fix any $a\in \R$ and $b\in (a, a+2/\sqrt N]$.
Let 
$$I_j=\{ \exists t\in [(j-1)/N, j/N]: a\leq W_t\leq b)\},$$
and 
$$A_k=\sum_{j=1}^k \bone_{I_j}.$$

\proclaim Lemma 3.1. There exist $c_1$ and $ c_2$ such that for all $\lam>0$,
$$\P(A_k\geq \lam \sqrt k )\leq c_1e^{-c_2\lam}.$$

\proof There is probability $c_3>0$ independent of $x$ such that 
$$\P^x(\sup_{s\leq 1/N} |W_s-W_0|
<1/\sqrt N)>c_3.$$ 
So by the strong Markov property applied at the first
$t\in [(j-1)/N, j/N]$ such that $ a\leq W_t\leq b$,
$$c_3\P^x(I_j)\leq \P^x(W_{j/N}\in [a-(1/\sqrt N), a+(3/\sqrt N)]).$$
This and the standard bound
$$\P^x(W_t\in [c,d])=\int_c^d \frac{1}{\sqrt{2\pi t}}e^{-|y-x|^2/2t}dy\leq \frac{1}{\sqrt {2\pi t}}|d-c|,$$
imply that
$$\P^x(I_j)\leq c_4\frac{1}{\sqrt N}\frac{1}{\sqrt{j/N}}=\frac{c_4}{\sqrt j}.$$
Therefore
$$\E^x A_k=\sum_{j=1}^k \P(I_j)\leq c_5 \sqrt k.\eqno(3.1)$$
By the Markov property, 
$$\E[A_k-A_i\mid \F_{i/n}]\leq 1+\E^{W(i/n)} A_k\leq c_6\sqrt k.\eqno(3.2)$$
Corollary I.6.12 of Bass (1995) can be applied to the sequence
$A_k/(c_7\sqrt{k})$, in view of (3.1) and (3.2). That 
result say that $\E \exp (c_8 \sup_k A_k/(c_7\sqrt{k})) \leq 2$
for some $c_8 >0$. This easily implies our lemma.
\qed

Fix an integer $N>0$. 
Let $R_{\ell m}= R_{\ell m}(N)$ be the rectangle defined by
$$R_{\ell m}= [\ell/N, (\ell+1)/N]\times [m/N^\al, (m+1)/N^\al], 
\qq 0\leq \ell\leq N, 
\quad -N^\al-1\leq m\leq N^\al.$$
Let $K$ be 
such that $N/K$ is an integer and $\sqrt N < N/K\leq \sqrt N+1$. 
Set
$$Q_{ik}=Q_{ik}(N) = [iK/N, (i+1)K/N]\times [k(K/N)^\al,(k+1)(K/N)^\al],$$
for $0\leq i\leq K$ and $ -(N/K)^\al-1 \leq k\leq (N/K)^\al$.
Note that $Q_{ik}(N) = R_{ik}(N/K) $ but it will
be convenient to use both notations.

\proclaim Proposition 3.2. 
Let $\al \in (1/2,1]$ and $\eps\in (0,1/16)$.
There exist $c_1, c_2$, and $c_3$ such that:
\item{(i)} there exists a set $D_N$ with $\P(D_N)\leq c_1N\exp(-c_2N^{\eps/2})$;
\item{(ii)} if $\omega\notin D_N$ and $f\in \sS_\al$, 
then there are at most $c_3 N^{(3/4)+(\eps/2)}$ 
rectangles $R_{\ell m}$ in $[0,1]\times [-1,1]$ which contain 
both a point of the graph of $f$ and a point of the graph of $W_t(\omega)$.

\proof
Let 
$$I_{ikj}=\{ \exists t\in [iK/N+(j-1)/N, iK/N+j/N]: 
k(K/N)^\al\leq W_t\leq (k+1)(K/N)^\al \},$$
$$A_{ik}=\sum_{j=1}^K \bone_{I_{ikj}},$$
and
$$C_{ik}=C_{ik}(N)=\{A_{ik}\geq K^{(1/2)+\eps}\}.$$
By Lemma 3.1 with $k=[K]$ and $\lam = K^\eps$, and the Markov property
applied at $kK/N$
we have $\P(C_{ik})\leq c_4\exp(-c_5 K^{\eps})$.

There are at most $c_6N^{(1/2)+(\al/2)}$ rectangles $Q_{ik}$, so if
$D_N=\cup_{i,k} C_{ik}$, where 
$0\leq i\leq K$ and $ -(N/K)^\al-1 \leq k\leq (N/K)^\al$,
then
$$\P(D_N)\leq c_7 N^{(1+\al)/2} \exp(-c_5 K^{\eps})
\leq c_7 N \exp(-c_8 N^{\eps/2}).$$

Now suppose $\omega\notin  D_N$. Let $f$ be
any function in $\sS_\al$. If $f$ intersects $Q_{ik}$ for some $i$ and
$k$, then $f$ might intersect $Q_{i,k-1}$ and $Q_{i,k+1}$.
But because $f\in\sS_\al$, it cannot
intersect $Q_{ir}$ for any $r$ such that $|r-k|>1$. Therefore
$f$ can intersect at most $3(K+1)$ of the $Q_{ik}$. 

Look at any one of the $Q_{ik}$ that $f$ intersects. Since $\omega\notin
D_N$, then there are at most $K^{(1/2)+\eps}$ integers $j$ that
are less than $K$ and for which the path of $W_t(\omega)$ intersects
$([iK/N+(j-1)/N,iK/N+j/N]\times[-1,1])\cap Q_{ik}$.
If $f$ intersects a rectangle $R_{\ell m}$,
then it can intersect a rectangle $R_{\ell r}$
only if $|r-m|\leq 1$, since $f\in\sS_\al$.
Therefore there are at most $3K^{(1/2)+\eps}$ rectangles
$R_{\ell m}$ contained
in $Q_{ik}$ which contain both a point of the graph of $f$
and a point of the graph of $W_t(\omega)$.

Since there are at most $3(K+1)$ rectangles $Q_{ik}$ which contain
a point of the graph of $f$, there are therefore at most
$$3(K+1) 3K^{(1/2)+\eps}\leq c_9 N^{(3/4)+(\eps/2)}$$ 
rectangles $R_{\ell m}$ that
contain both a point of the graph of $f$ and a point of the graph of $W_t(\omega)$.
\qed

\bs

We can now iterate this to obtain a better estimate.

\proclaim Proposition 3.3. Fix $\al\in(1/2,1]$
and $\delta,\eta>0$. There exist $c_1$ and $N_0$
such that if $N\geq N_0$:
\item{(i)} there exists a set $E$ with $\P(E)\leq \eta$;
\item{(ii)} if $\omega\notin E$ and  $f\in \sS_\al$,
then there are at most
$c_1 N^{(1/2)+\delta}$ rectangles $R_{\ell m}(N)$ contained in $[0,1]\times [-1,1]$ 
which contain
both a point of the graph of $f$ and a point of the graph of $W_t(\omega)$.

\proof
For any $\eps$, the quantity $c_1N\exp(-c_2N^{\eps/2})$ is summable.
First choose  $\eps\in (0,\delta/4)$ and
then choose $N_1$ large so that, using Proposition 3.2 and its notation,
$$\sum_{N=N_1}^\infty \P(D_N) \leq \sum_{N=N_1}^\infty
c_1N\exp(-c_2N^{\eps/2}) <\eta.$$
Let $E=\cup_{N=N_1}^\infty D_N$.

Fix $\omega\notin E$.
Suppose $N$ is large
enough so that $\sqrt N\geq 2N_1$. 
Recall the definition of $K$ and note that $N/K$
differs from $\sqrt{N}$ by at most $1$.
Then by Proposition 3.2 applied with $N/K$,
there are at most  $c_2(\sqrt N)^{(3/4)+\eps}$ rectangles 
$R_{ik}(N/K)$
that contain both  a point of the graph of $f$ and a point of the graph of $W_t(\omega)$. 
Recall the definitions of the events $C_{ik}$ and $D_N$ from
Proposition 3.2 and its proof. Since we are assuming that
$\omega \notin E$, we also have $\omega \notin C_{ik}(N)$
for any $i,k$. This implies that
inside each rectangle $R_{ik}(N/K)$, there are at most
$c_3(\sqrt N)^{(1/2)+\eps}$ rectangles $R_{\ell m}(N)$ that contain both
a point of the graph of $f$ and a point of the graph of $W_t(\omega)$.
Thus there are at most 
$$c_4(\sqrt N)^{(3/4)+\eps}(\sqrt N)^{(1/2)+\eps}=c_4N^{(5/8)+\eps}$$
rectangles $R_{\ell m}(N)$ that contain both
a point of the graph of $f$ and a point of the graph of $W_t(\omega)$.
 
We continue
iterating: take $N$ large so that $N\geq (4N_1)^4$. There are 
$c_4(\sqrt N)^{(5/8)+\eps}$ rectangles $R_{\ell m}(N/K)$
that contain both
a point of the graph of $f$ and a point of the graph of $W_t(\omega)$.
Each of these contains at most $c_5(\sqrt N)^{(1/2)+\eps}$ 
rectangles $R_{\ell m}(N)$ that contain both
 a point of the graph of $f$ and a point of the graph of $W_t(\omega)$,
for a total of 
$$c_6(\sqrt N)^{(5/8)+\eps}(\sqrt N)^{(1/2)+\eps}=c_6N^{(9/16)+\eps}$$
rectangles $R_{\ell m}(N)$. 

Continuing, if $N$ is large enough, we can get the exponent of 
$N$ as close to $(1/2)+\eps$ as we like. In particular, by a finite
number of iterations, we can get the exponent less than $(1/2)+\delta$.
\qed

\bs

Recall the definition of $p(t,x,y)$ in (2.1).

\proclaim Lemma 3.4. If $\norm{f-g}\leq \eps$, then for some constant
$c_1$ and all $\eps < \frac12$,
$$\int_0^1 |p(t,0,f(t))-p(t,0,g(t))|\, dt\leq c_1\eps\log (1/\eps).$$

\proof For $t\leq \eps^2$, we use the estimate $p(t,0,x)\leq c_2 t^{-\half}$
and obtain
$$\int_0^{\eps^2} |p(t,0,f(t))-p(t,0,g(t))|\, dt
\leq 2c_2\int_0^{\eps^2} \frac{1}{\sqrt t}dt\leq c_3\eps.$$
For $t\geq \eps^2$, note that
$$\Bigl|\frac{\del p(t,0,x)}{\del x}\Bigr|
=c_4 t^{-\half} \frac{|x|}{t}e^{-x^2/2t}= c_4 t^{-1}
\frac{|x|}{\sqrt t} e^{-x^2/2t}\leq c_5 t^{-1},$$
since $|y|e^{-y^2/2}$ is bounded. We then obtain
$$\int_{\eps^2}^1  |p(t,0,f(t))-p(t,0,g(t))|\, dt
\leq \int_{\eps^2}^1 |f(t)-g(t)| c_5 t^{-1}dt
\leq c_5\eps \int_{\eps^2}^1  t^{-1}\,dt=c_6\eps\log(1/\eps).$$
Adding the two integrals proves the lemma.
\qed

\proclaim Proposition 3.5. Let $f$ and $g$  be two functions
with 
$$\sup_{(j-1)/N\leq t\leq j/N} |f(t)-g(t)|\leq \delta.$$
Then, for all $\lam >0$, 
$$\P\big(|(L^f_{j/N}-L^f_{(j-1)/N}) - (L^g_{j/N}-L^g_{(j-1)/N})|
\geq \lam N^{-1/4}(\delta\log (1/\delta))^\half\big)\leq c_1 e^{-c_2\lam}.$$

\proof Write $s$ for $(j-1)/N$ and $A^f_t=L^f_{s+t}-L^f_s$,
$A^g_t=L^g_{s+t}-L^g_s$. We have
for $s\leq r\leq t\leq s+(1/N)$,
$$\E[A^f_t-A^f_r\mid \F_r]=\E^{W_r} A^f_{t-r}
\leq \sup_z \E^z A^f_{1/N}.$$
But for any $z$,
$$\E^z A^f_{1/N}=\int_0^{1/N} p(t,z,f(t))\, dt
\leq \int_0^{1/N} \frac{1}{\sqrt t}dt\leq c_3 N^{-\half}.$$
We have a similar bound for $\E^z A^g_{1/N}$.
For the difference, we have
$$|\E[ (A^f_t-A^g_t)-(A^f_r -A^g_r)\mid \F_r]|
=|\E^{W_r} [A^f_{t-r}-A^g_{t-r}]|.$$
However, for any $z$,
$$\eqalign{|\E^z[[A^f_{t-r}-A^g_{t-r}]| &=\Bigl|\int_s^{s+t-r} [p(u,z,f(u))
-p(u, z,g(u))] du\Bigr|\cr
&\leq \int_0^1 |p(u, 0, \wt f(u))-p(u,0,  \wt g(u))|\, du,\cr}$$
where we define $\wt f(u)=f(u)-z$ for all $u$ and we
define $\wt g(u)=g(u)-z$ if $s\leq u\leq s+(t-r)$ and
$\wt g(u)=\wt f(u)$ otherwise. So $\norm{\wt f(u)-\wt g(u)}\leq \delta$,
and by Lemma 3.4, 
$$|\E^z[[A^f_{t-r}-A^g_{t-r}]| \leq c_4\delta\log(1/\delta).$$

Our result now follows by Lemma 2.1. 
\qed

\proclaim Theorem 3.6. For any $\al \in (1/2,1]$,
there exists $\wt L^f_t$ such that 
\item{(i)} for each $f\in \sS_\al$, we have 
$\wt L^f_t=L^f_t$ for all $t$, a.s., 
\item{(ii)} with probability one, $ f\to \wt L^f_1$ is a continuous
map on $\sS_\al$ with respect to the supremum norm, and
\item{(iii)} with probability one, $\sup_{f\in \sS_\al} \wt L^f_1 <\infty$.

\proof

\noindent{\it Step 1}. In this step, we will define
and analyze a countable dense family of functions in $\sS_\al$.

Let $N=2^n$ and let $T_n$ denote the class of functions $f$ in $\sS_\al$ 
such that
on each interval $[(j-1)/N, j/N]$ the function $f$ is linear with
slope either $N^{1-\al}$ or $-N^{1-\al}$ and $f(j/N)$ is a multiple 
of $1/N^\al$ for each $j$.
Note that the collection of all functions which are piecewise linear
with these slopes contains some functions which are not in $\sS_\al$-- such
functions do not belong to $T_n$.

Consider any element $h$ of $\sS_\al$. Let $h^{(n)}$ denote
a function in $T_n$ which approximates $h$ in the following sense.
We will define $h^{(n)}$ inductively on intervals of the form
$[(j-1)/N, j/N]$. 
First we take the initial value $h^{(n)}(0)$ to be the
closest integer multiple of $1/N^\al$ to $h(0)$ (we take the smaller
value in case of a tie).
The slope of $h^{(n)}$ is chosen to be positive on
$[0, 1/N]$ if and only if $h^{(n)}(0) \leq h(0)$. Once the function
$h^{(n)}$ has been defined on all intervals $[(j-1)/N, j/N]$,
$j = 1,2,\dots, k$, we choose the slope of $h^{(n)}$
on $[k/N, (k+1)/N]$ to be $N^{1-\al}$ if and only if
$h^{(n)}(k/N) \leq h(k/N)$. Strictly speaking, our definition
generates some functions with values in $[-1-1/N^\al, 1+ 1/N^\al]$
rather than in $[-1,1]$ and so $h^{(n)}$ might not belong to
$\sS_\al$. We leave it to the reader to check that
this does not affect our arguments.

We will argue that $|h^{(n)}(t) - h(t)| \leq 2/N^\al$
for all $t$. This is true for $t=0$ by definition.
Suppose that $1/N^\al \leq |h^{(n)}(t) - h(t)| \leq 2/N^\al$ for some $t=j/N$.
Then the fact that both functions belong to $\sS_\al$
and our choice for the slope of $h^{(n)}$ easily
imply that the absolute value of the difference between the two functions
will not be greater at time $t= (j+1)/N$ than at time $t=j/N$.
An equally elementary argument shows that in the case
when $ |h^{(n)}(t) - h(t)| \leq 1/N^\al$, the distance
between the two functions may sometimes increase but will
never exceed $2/N^\al$. The induction thus proves the claim
for all times $t$ of the form $t=j/N$. An extension to all
other times $t$ is easy.

Later in the proof we will need to consider the difference
between $h^{(n)}$ and $h^{(n+1)}$. 
First let us restrict our attention to the interval
$[\ell/N, (\ell+1)/N]$. The estimates from the previous paragraph
show that $|h^{(n)}(t) - h^{(n+1)}(t)| \leq 4/N^\al$
on this interval. 
Let 
$$F_{h,\ell }=\{|(L^{h^{(n)}}_{(\ell+1)/N}-L^{h^{(n)}}_{\ell/N})
-(L^{h^{(n+1)}}_{(\ell+1)/N}-L^{h^{(n+1)}}_{\ell/N})|
\geq N^{-(1/4)-(\al/2)+\eps}\}.$$
By Proposition 3.5 with $\lam = N^{\eps}$, for any $h\in \sS_\al$, $\ell$
and $n$,
$$\P(F_{h,\ell })\leq c_1 \exp(-c_2N^{\eps}).$$

There are only $N+1$ integers $\ell$ with $0 \leq \ell\leq N$.
For a fixed $\ell$, there are no more than $3N^\al$ possible values 
of $h^{(n)}(\ell/N)$, and the same is true for $h^{(n)}((\ell+1)/N)$.
The analogous
upper bound for the number of possible values for each of
$h^{(n+1)}(\ell/N)$, $h^{(n+1)}((\ell+1/2)/N)$ and
$h^{(n+1)}((\ell+1)/N)$ is 
$6N^\al$. Hence, if we let
$$G_N=\bigcup_{h\in \sS_\al} \bigcup_{0 \leq \ell\leq N} F_{h,\ell },$$
then  
$$\P(G_N)\leq  c_3 N^6  \exp(-c_2N^{\eps}).$$

We will derive a similar estimate for $f^{(n)}$ and $h^{(n)}$,
where $f,h \in \sS_\al$. Let us assume that $\| f - h \| \leq 1/N^\al$.
Then $|f^{(n)}(t) - h^{(n)}(t)| \leq 5/N^\al$ for all $t$.
If we define
$$\wt F_{f,h,\ell }=\{|(L^{f^{(n)}}_{(\ell+1)/N}-L^{f^{(n)}}_{\ell/N})
-(L^{h^{(n)}}_{(\ell+1)/N}-L^{h^{(n)}}_{\ell/N})|
\geq N^{-(1/4)-(\al/2)+\eps}\}.$$
then
$$\P(\wt F_{f,h,\ell })\leq c_7 \exp(-c_8N^{\eps}).$$
Next we let
$$\wt G_N=\bigcup_{f,h\in \sS_\al} 
\bigcup_{0 \leq \ell\leq N} \wt F_{f,h,\ell }.$$
Counting all possible paths $f^{(n)}$ and $h^{(n)}$
yields an estimate analogous to the one for $G_N$,
$$\P(\wt G_N)\leq  c_9 N^5  \exp(-c_8 N^{\eps}).$$

\medskip
\noindent{\it Step 2}.
In this step, we will prove uniform continuity of $f \to L^f_1$
on the set $T_\infty=\bigcup_{n=1}^\infty T_n$.

Fix arbitrarily small $\eta,\beta>0$. Choose $\eps>0$
so small that
$(1/4)-(\al/2)+2\eps < 0$.
Recall the events $D_N$ from Proposition 3.2. Since
$\sum_N (\P(D_N)+\P(G_N)+\P(\wt G_N)) < \infty$,
we can take $N_0$ sufficiently large so that $\P(H)\leq \eta$, where
$H=\bigcup_{N=N_0}^\infty (D_N\cup G_N\cup \wt G_N)$. 
Without loss of generality we may take $N_0$ to be an integer power of
$2$, say $N_0=2^{n_0}$.

Fix an $\omega\notin H$. Consider any $f,h\in T_\infty$
with $\| f - h\| \leq 1/N_0^\al$.
Note that
$$|L^h_1-L^{h^{(n_0)}}_1|\leq \sum_{n=n_0}^\infty
|L^{h^{(n+1)}}_1-L^{h^{(n)}}_1|,\eqno(3.3)$$
and
$$|L^{h^{(n+1)}}_1-L^{h^{(n)}}_1|\leq \sum_{m=1}^{2^n} 
|(L^{h^{(n+1)}}_{(m+1)/2^n}-L^{h^{(n+1)}}_{m/2^n}) -
(L^{h^{(n)}}_{(m+1)/2^n}-L^{h^{(n)}}_{m/2^n})|. \eqno (3.4)$$

Consider $2^n=N\geq N_0$. Since $\omega \notin \bigcup _{N\geq N_0} D_N$,
Proposition 3.3 implies that
there are at most $c_1 N^{(1/2)+\eps}$ values of $m$
for which there is a rectangle $R_{mi}$ in which there is a point of the 
graph of $h^{(n)}$ or of $h^{(n+1)}$  and a point of the graph of $W_t(\omega)$. 
So there are no more than $c_1N^{(1/2)+\eps}$ summands
on the right hand side of (3.4) that are non-zero. 

For a value of $m$ for which the summand on the right
hand side is nonzero, it is at most $N^{-(1/4)-(\al/2)+\eps}$,
because $\omega\notin \bigcup _{N\geq N_0} G_N$.
Multiplying the number of nonzero summands by the the largest value
each summand can be, we obtain
$$\eqalignno{|L^{h^{(n+1)}}_1-L^{h^{(n)}}_1|&\leq c_1N^{(1/2)+\eps}N^{-(1/4)-(\al/2)+\eps}\cr
&=c_1 N^{(1/4)-(\al/2)+2\eps}= c_1 (2^n)^{(1/4)-(\al/2)+2\eps}.& (3.5)\cr}$$
We have assumed that $\eps$ is so small that
$(1/4)-(\al/2)+2\eps <0$, so the bound in (3.5) is summable in $n$.
We increase $n_0$, if necessary, so that
$\sum_{n\geq n_0} c_1 (2^n)^{(1/4)-(\al/2)+2\eps} \leq \beta/3$.
Then (3.3) implies that
$$|L^h_1-L^{h^{(n_0)}}_1|\leq \beta/3.$$
Similarly,
$$|L^f_1-L^{f^{(n_0)}}_1|\leq \beta/3.$$

A similar reasoning will give us a bound for
$|L^{f^{(n_0)}}_1 - L^{h^{(n_0)}}_1|$.
We have
$$|L^{f^{(n_0)}}_1 - L^{h^{(n_0)}}_1|
\leq \sum_{\ell=1}^{2^n}
|(L^{f^{(n)}}_{(\ell+1)/N}-L^{f^{(n)}}_{\ell/N})
-(L^{h^{(n)}}_{(\ell+1)/N}-L^{h^{(n)}}_{\ell/N})|.$$
First, the number of non-zero summands is bounded
by $c_1N_0^{(1/2)+\eps}$, for the same reason as above.
We have assumed that $\| f - h\| \leq 1/N_0^\al$,
so, in view of the fact that $\omega \notin \bigcup _{N\geq N_0}\wt G_N$,
the size of a non-zero summand is bounded by
$N_0^{-(1/4)-(\al/2)+\eps}$.
Hence,
$$|L^{f^{(n_0)}}_1 - L^{h^{(n_0)}}_1|
\leq  c_1N_0^{(1/2)+\eps}N_0^{-(1/4)-(\al/2)+\eps}
= c_1 (2^{n_0})^{(1/4)-(\al/2)+2\eps}\leq \beta/3.$$

By the triangle inequality, with probability greater than $1-\eta$,
$$|L^f_1-L^h_1| \leq \beta$$
if $f,h\in T_\infty$ and $\| f - h\| \leq 1/N_0^\al \df \delta(\beta)$.
We now fix
an arbitrarily small $\eta_0>0$ and a sequence $\beta_k \to 0$,
and find $\delta(\beta_k)>0$ such that 
with probability greater than $1-\eta_0/2^k$,
$$|L^f_1-L^h_1| \leq \beta_k,$$
if $f,h\in T_\infty$ and $\| f - h\| \leq \delta(\beta_k)$.
This implies that, with probability greater than $1-\eta_0$, the function
$f \to L^f_1$ is uniformly continuous on $T_\infty$.
Since $\eta_0$ is arbitrarily small, the uniform continuity
is in fact an almost sure property, although the modulus
of continuity may depend on $\omega$.

For an arbitrary $f\in\sS_\al$, define
$\wt L^f=\lim_{n\to\infty} L^{f^{(n)}}_1$. 
By Corollary 2.5, $L^f=\wt L^f$ a.s. Therefore
$\wt L^f$ is a version of $L^f$. 

Since the function $f\to L^f_1$ is uniformly continuous
on $T_\infty$, its extension to $\sS_\al$ is uniformly
continuous with the same (random) modulus of continuity.
The family $\sS_\al$ is
equicontinuous, hence a compact set with respect to $\norm{\cdot}$.
Therefore the supremum  of $\wt L^f_1$ over $\sS_\al$ is finite, a.s.
\qed

\noindent{\bf Remark 3.7}. It is rather easy to see
that, with probability one, $f\to \wt L^f_t$ is actually jointly
continuous on $\sS\times [0,1]$.
To see this, note that in the proof of Proposition 3.5 we used Proposition 2.1,
so what we actually proved was that
$$\P\left(\sup_{(j-1)/n\leq t\leq j/n}
|(L^f_t-L^f_{(j-1)/n})-(L^g_t-L^g_{(j-1)/n})|
\geq \lam N^{-1/4}(\delta \log(1/\delta))^\half\right)
\leq e^{-c_1\lam}. $$
If we replace (3.4) by
$$\sup_t|L^{h^{(n+1)}}_t-L^{h^{(n)}}_t|\leq
\sum_{m=1}^{2^n}\sup_{m/2^n\leq t\leq (m+1)/2^n} 
|(L^{h^{(n+1)}}_t-L^{h^{(n+1)}}_{m/2^n}) -
(L^{h^{(n)}}_t-L^{h^{(n)}}_{m/2^n})|,$$
then proceeding as in the proof of Theorem 3.6,
we obtain the joint continuity.
\bs

We will show that, in a sense, $\sup_{f\in \sS_{\al}}L^f_1=\infty$,
a.s., if $\al < \half$. This statement is quite intuitive
-- one would like to let $f(\omega)=W_t(\omega)$ so that $L^f_1(\omega)=\infty$ -- 
but we have not defined the local time simultaneously
for all $f\in \sS_{\al}$, and there is a difficulty with the
number of null sets. Theorem 3.6 suggests that
the question of joint existence is tied to the question
of the finiteness of the supremum, so we have
to express our result in a different way.

\proclaim Theorem 3.8. Suppose $\al<\half$. 
Then there exists a countable family $F\subset \sS_\al$ such that
$\sup_{f\in F}L^f_1=\infty$ a.s.

\proof 
Let $\ell^x_t$ be the ordinary local time at $x$ for Brownian motion.
It is well known (see Karatzas and Shreve (1994)) that
there exists a version of this process which is jointly
continuous in $x$ and $t$.

Suppose that a piecewise linear function $f$ is equal to $y$ on an interval
$[s,t]$. Then Proposition 2.2 and a similar well known
result for $\ell^y$ show that with probability one,
for all $u\in [s,t]$,
$$L^f_u - L^f_s = \ell^y_u - \ell^y_s.$$

Fix $\al \in (0,1/2)$.
Let $F$ be the countable family of all functions $f$
defined on the interval $[0,1]$ such that for some
integers $n=n(f)$ and $m=m(f)$,
on each interval of the form
$[(j-1)/n, (j-\frac12)/n]$ the function $f$ is a 
constant multiple of $2^{-m}$,
$f$ is linear on the intervals $[(j-\frac12)/n, j/n]$, and $f\in \sS_\al$.
Then, with probability one, for all $j$, all $f\in F$ and $n=n(f)$,
$$ L^{f((j-1)/n)}_{(j-(1/2))/n}- L^{f((j-1)/n)}_{(j-1))/n}
=\ell^{f((j-1)/n)}_{(j-(1/2))/n}- \ell^{f((j-1)/n)}_{(j-1))/n}. \eqno(3.6)$$
In the rest of the proof we assume that this assertion
and the joint continuity of $\ell^x_t$ hold for all $\omega$.

Let 
$$T=\inf\{t: |W_t|\geq 1\hbox{ or } \exists r,s\leq t
\hbox{ such that } |W_r-W_s|\geq (\textstyle{\frac14} |r-s|)^\al\}.\eqno(3.7)$$
By the well-known results on the modulus
of continuity for Brownian motion, $T>0$ a.s.

Let $\eps>0$. There exists $\delta$ such that $\P(T<\delta)<\eps$.
Fix $n$.
On the interval $[(j-1)/n,(j-\frac12)/n]$,
let $ f_1(t)= W((j-1)/n)$. On
the interval 
$[(j-\frac12)/n,j/n]$ let $ f_1(t)$ be linear with
$ f_1(j/n)=W(j/n)$. Let $f_2(t)= f_1(t)$ for $t\leq \delta/2$ and
constant for $t\geq \delta/2$.

It is quite easy to show that $f_2\in \sS_\al$ for each $\omega$ in  
the set $\{T>\delta\}$ using the definition (3.6) of $T$.
By the Markov property, the random variables
$$X_j=\ell^{f_2((j-1)/n)}_{(j-(1/2))/n}- \ell^{f_2((j-1)/n)}_{(j-1))/n}$$
form an independent sequence, and by Brownian scaling,
$Y_j=\sqrt{2n} X_j$ has the same distribution as $\ell^0_1$. Let $c_1=\E \ell^0_1$.
By Chebyshev's inequality,
$$\P\Big(\Bigl| \sum_{j=1}^{[\delta n/2]} 
(Y_j-c_1)\Bigl|\geq c_1\delta n/4\Big)
\leq \frac {[\delta n/2] {\rm Var}\, Y_1}{ (c_1\delta n/4)^2}
\leq \frac{c_2\E(\ell^0_1)^2}{\delta n}=\frac{c_3}{\delta n}.$$
Take $n$ large so that $c_3/(\delta n)<\eps$. Then there exists a set $A_n$
of probability at most $2\eps$ such that if $\omega\notin A_n$, then
$T(\omega)\geq \delta$ and 
$$\sum_{j=1}^{[\delta n/2]} X_j\geq c_4\sqrt {\delta n}.$$

We now choose $m$ large and 
find $f_3\in F$ so that on each interval
$[(j-1)/n, (j-\frac12)/n]$ the function $f_3$ is a multiple of $2^{-m}$,
$f_3$ is linear on the intervals $[(j-\frac12)/n, j/n]$, and
$$\sum_{j=1}^{[\delta n/2]} \Big[ \ell^{f_3((j-1)/n)}_{(j-(1/2))/n}- 
\ell^{f_3((j-1)/n)}_{(j-1))/n}\Big]\geq c_4\sqrt {\delta n}/2;$$
this is possible by the joint continuity of $\ell^x_t$.

By (3.6) we can replace $\ell$ by $L$
in the last formula, so
$$L^{f_3}_1 \geq \sum_{j=1}^{[\delta n/2]} \Big[ L^{f_3}_{(j-(1/2))/n}- 
L^{f_3}_{(j-1))/n}\Big]\geq c_4\sqrt {\delta n}/2.$$
We conclude that
$$\sup_{f \in F} L^f_1 \geq c_4\sqrt {\delta n}/2,$$
with probability greater than or equal to $1-2\eps$.
Since $n$ and $\eps$ are arbitrary, the proposition is proved.
\qed

\centerline{\bf References}
\medskip

\item{1.}R.F. Bass (1995).  {\sl Probabilistic Techniques in Analysis}. 
Springer-Verlag, New York.

\item{2.}R.F. Bass and D. Khoshnevisan (1992). Local times on curves and uniform 
invariance principles. {\sl Probab. Theory
Related Fields  \bf 92}, 465--492.

\item{3.}R.F. Bass and K. Burdzy (1999). Stochastic bifurcation models. {\sl Ann. Probab.  \bf 27}, 50--108. 

\item{4.}K. Burdzy and  J. San Mart\'in (1995). Iterated law of iterated logarithm. {\sl Ann. Probab.  \bf 23},
1627--1643.

\item{5.}G.F. Clements (1963). Entropies of several sets of real valued functions. {\sl Pacific J. Math.  \bf 13},   1085--1095. 

\item{6.}B.  Davis, (1998). Distribution of Brownian local time on curves. {\sl Bull. London Math. Soc.  \bf 30}, 182--184.

\item{7.} I.~Karatzas and S.~Shreve (1994) {\it Brownian Motion and
Stochastic Calculus, Second Edition}. Springer-Verlag,
New York.

\item{8.}M. Ledoux and M.  Talagrand (1991). 
{\sl Probability in Banach spaces. Iso\-per\-i\-me\-try and Processes.} 
Springer-Verlag, Berlin.

\item{9.}M.B. Marcus and J.  Rosen (1996). Gaussian chaos and sample path properties of additive functionals of symmetric
Markov processes. {\sl Ann. Probab. \bf  24}, 1130--1177.

\item{10.}D. Revuz and M. Yor (1994). {\sl Continuous Martingales
and Brownian Motion, 2nd ed.} Springer-Verlag, Berlin.

\vskip1truecm
 
\parskip=0pt
\hbox{\vbox{
\obeylines
Richard F.~Bass
Department of Mathematics
University of Connecticut
Storrs, CT 06269
e-mail: bass@math.uconn.edu
{\ }
}
\hskip-8truecm
\vbox{
\obeylines
Krzysztof Burdzy
Department of Mathematics
University of Washington
Box 354350
Seattle, WA 98195-4350
e-mail: burdzy@math.washington.edu}
}

\bye

%% file: rbmacros.tex
\magnification=\magstep1

\baselineskip=14pt

\expandafter\ifx\csname bookmacros.tex\endcsname\relax \else \fi
%
\expandafter\edef\csname bookmacros.tex\endcsname{%
       \catcode`\noexpand\@=\the\catcode`\@\space}
\catcode`\@=11

\def\undefine#1{\let#1\undefined}
\def\newsymbol#1#2#3#4#5{\let\next@\relax
 \ifnum#2=\@ne\let\next@\msafam@\else
 \ifnum#2=\tw@\let\next@\msbfam@\fi\fi
 \mathchardef#1="#3\next@#4#5}
\def\mathhexbox@#1#2#3{\relax
 \ifmmode\mathpalette{}{\m@th\mathchar"#1#2#3}%
 \else\leavevmode\hbox{$\m@th\mathchar"#1#2#3$}\fi}
\def\hexnumber@#1{\ifcase#1 0\or 1\or 2\or 3\or 4\or 5\or 6\or 7\or 8\or
 9\or A\or B\or C\or D\or E\or F\fi}

\font\tenmsa=msam10
\font\sevenmsa=msam7
\font\fivemsa=msam5
\newfam\msafam
\textfont\msafam=\tenmsa
\scriptfont\msafam=\sevenmsa
\scriptscriptfont\msafam=\fivemsa
\edef\msafam@{\hexnumber@\msafam}
\mathchardef\dabar@"0\msafam@39

\font\tenmsb=msbm10
\font\sevenmsb=msbm7
\font\fivemsb=msbm5
\newfam\msbfam
\textfont\msbfam=\tenmsb
\scriptfont\msbfam=\sevenmsb
\scriptscriptfont\msbfam=\fivemsb
\edef\msbfam@{\hexnumber@\msbfam}
\def\Bbb#1{{\fam\msbfam\relax#1}}
\def\widehat#1{\setbox\z@\hbox{$\m@th#1$}%
 \ifdim\wd\z@>\tw@ em\mathaccent"0\msbfam@5B{#1}%
 \else\mathaccent"0362{#1}\fi}
\def\widetilde#1{\setbox\z@\hbox{$\m@th#1$}%
 \ifdim\wd\z@>\tw@ em\mathaccent"0\msbfam@5D{#1}%
 \else\mathaccent"0365{#1}\fi}
\font\teneufm=eufm10
\font\seveneufm=eufm7
\font\fiveeufm=eufm5
\newfam\eufmfam
\textfont\eufmfam=\teneufm
\scriptfont\eufmfam=\seveneufm
\scriptscriptfont\eufmfam=\fiveeufm

\newsymbol\centerdot 1205
\newsymbol\square 1003

\csname bookmacros.tex\endcsname

\def\R{{\Bbb R}}
\def\E{{\Bbb E}}
\def\P{{\Bbb P}}

\def\F{{\cal F}}

\def\lam{{\lambda}}

\def\al{{\alpha}}

\def\proof{{\medskip\noindent {\bf Proof. }}}

\def\qed{{\hfill $\square$ \bigskip}}
\def\subsec#1{{\bigskip\noindent \bf{#1}.}}

\def\section#1#2{{\bigskip\bigskip \centerline{\bf #1. #2}\bigskip}}

\def\chapter#1#2{{\bigskip\bigskip \centerline{#1. #2}\bigskip}}
\def\cite#1{{[#1]}}

\def\eps{\varepsilon}

\def\angel#1{{\langle#1\rangle}}

\def\norm#1{\Vert #1 \Vert}

 \def\qq {\qquad}
\def\frac#1#2{{#1\over #2}}
\def\del{{\partial}}
\def\wt{\widetilde}

\def\ni{\noindent }
\def\ms{\medskip}
\def\bs{\bigskip}
\def\cl#1{\centerline{#1}}

\def \half {{{1/ 2}}}

\parindent=30pt

\def\ftN{\rl\kern-0.13em\rN}

\def\sgn{{\mathop {{\rm sgn\, }}}}

\def\sqr#1#2{{\vcenter{\vbox{\hrule height.#2pt
        \hbox{\vrule width.#2pt height#1pt \kern#1pt
           \vrule width.#2pt}
        \hrule height.#2pt}}}}